\documentclass[12pt]{amsart}
\usepackage{amsfonts}

\theoremstyle{definition}

\theoremstyle{remark}

\newcommand{\const}{\mathop{\rm const}\limits}

\newcommand{\grad}{\mathop{\rm grad}\limits}

\newcommand{\diam}{\mathop{\rm diam}\limits}

\newcommand{\supp}{\mathop{\rm supp}\limits}

\newcommand{\order}{\mathop{\rm order}\limits}

\begin{document}

\begin{center}

{\bf POINCAR\'E TYPE INEQUALITIES  FOR TWO DIFFERENT  \\

\vspace{4mm}

 BILATERAL GRAND LEBESGUE SPACES} \par

\vspace{3mm}

{\bf E. Ostrovsky}\\

e-mail: eugostrovsly@list.ru \\

\vspace{3mm}

{\bf L. Sirota}\\

e-mail: sirota3@bezeqint.net \\

\vspace{4mm}

 Abstract. \\

 \vspace{3mm}

 In this paper we obtain the non-asymptotic inequalities of Poincar\'e type
between function  and its weak  gradient belonging  the so-called
Bilateral Grand Lebesgue Spaces over general metric measurable space. We also prove
the sharpness of these inequalities. \\

\end{center}

\vspace{3mm}

2000 {\it Mathematics Subject Classification.} Primary 37B30,
33K55; Secondary 34A34, 65M20, 42B25.\\

\vspace{3mm}

{\it Key words and phrases:} Metric measure space, ball, diameter, exact constant,
 norm, Grand and ordinary Lebesgue Spaces and norms, integral and other linear operator, factorable estimation,
Poincar\'e-Lebesgue and Poincar\'e-Lipshitz spaces, sets and inequalities, exact estimations, module of continuity. \\

\vspace{3mm}

\section{Introduction}

\vspace{3mm}

 Let $  (X, M, \mu, d) $  be metric measurable space with finite non-trivial measure $ \mu:  \ 0 < \mu(X) < \infty $
and also with finite non-trivial distance function $ d = d(x,y): $

$$
 0 < \diam(X) := \sup_{x,y \in X} d(x,y) < \infty.
$$

  Define also for arbitrary numerical measurable function $ u: X \to R $ the following average

$$
u_X = \frac{1}{\mu(X)} \int_X u(x) \ d \mu(x),
$$

$$
||u||_p = \left[ \int_X |u(x)|^p \ d \mu(x)   \right]^{1/p}, \ p = \const \in [1, \infty],
$$
$ g(x) = \nabla[u](x)  $ will denote a so-called {\it minimal weak upper gradient} of the function $ u(\cdot), $ i.e.
the (measurable) minimal function such that for any rectifiable curve $ \gamma: [0,1] \to X $

$$
|u(\gamma(1)) - u(\gamma(0))| \le  \int_{\gamma} g(s) \ ds.
$$

 Note that if the function $ u(\cdot) $ satisfies the Lipschitz condition:

$$
|u(x) - u(y)| \le L \cdot d(x,y), \ 0 \le L = \const < \infty,
$$
then the function $ g(x) = \nabla[u](x)  $ there exists and is bounded:  $ g(x) \le L. $\par

\vspace{3mm}

 "The term Poincar\'e type inequality is used, somewhat loosely,
to describe a class of inequalities that generalize the classical Poincar\'e
inequality"

$$
\int_D |u(z)|^p \ dx \le A_m(p,D) \int_D | \ |\grad u(z)|^p \ | \ dz, \ A_m(p,D) = \const < \infty, \eqno(1.0)
$$
see \cite{Adams1}, chapter 8,  p.215, and the source work of Poincar\'e \cite{Poincare1}. \par
 A particular case  done by Wirtinger:  \hspace{1mm}
"an inequality ascribed to Wirtinger",  \ see \cite{Mitrinovich1}, p. 66-68;
see also \cite{Ostrovsky104}, \cite{Takahasi1}, \cite{Takahasi2}. \par

 In the inequality (1.0) $  D $ may be  for instance open bounded non empty convex subset of the whole
space $  R^m $  and has a Lipschitz or at last H\"older  boundary, or consists on the finite union of these
 domains, and $ |\grad u(z)|  $ is ordinary Euclidean $ R^m $ norm of "natural"  distributive gradient
of the differentiable a.e. function $  u. $ \par

\vspace{3mm}

 The {\it generalized} Sobolev's norm, more exactly, semi-norm $  ||f|| W_p^1 $ of a "weak differentiable" function $  f: \ f: X \to R $
may be defined by the formula

$$
||f|| W_p^1 \stackrel{def}{=} \left[  \int_X  |\nabla f|^p \ d \mu(x)  \right]^{1/p} = || \nabla f||_p.
$$

\vspace{3mm}

 We will call "the Poincar\'e inequality" , or more precisely  "the Poincar\'e $ (L(p), L(q))  $ inequality" more general
inequalities of the forms

$$
\mu(X)^{-1/q} ||u - u_X||_q \le K_P(p,q) \ \diam(X) \ \mu(X)^{-1/p} \ ||\nabla u||_p =
$$

$$
 K_P(p,q) \ \diam(X) \ \mu(X)^{-1/p} \ ||u||W_p^1, \eqno(1.1)
$$
where

$$
 1 \le p < s = \const > 1, \ 1 \le q <  \frac{ps}{s-p} \ -  \eqno(1.1a)
$$
the Poincar\'e-Lebesgue-Riesz version;
or in the case when $  p > s $ and after (possible) redefinition
of the function $ u = u(x) $ on a set of measure zero

$$
|u(x) - u(y)| \le K_L(s,p) \  d^{1 - s/p}(x,y) \ \mu(X) \ || \nabla u||_p =
$$

$$
 K_L(s,p) \ \mu(X) \ ||u||W_p^1 \  -  \eqno(1.2)
$$
the Poincar\'e-Lipshitz version; the case $ p = s $ in our setting of problem, indeed, in the terms of
Orlicz's spaces and norms, is considered in \cite{Koskela1}. \par

 The last inequality (1.2) may be reformulated in the terms of the module of continuity of the function $  u: $

$$
\omega(u, \tau) := \sup_{x,y: d(x,y) \le \tau} |u(x) - u(y)|, \ \tau \ge 0.
$$
Namely,

$$
\omega(u, \tau)  \le K_L(s,p) \  \tau^{1 - s/p} \ \mu(X) \  ||u||W_p^1||.
$$

\vspace{3mm}

 We will name following the authors of articles \cite{Heinonen1}, \cite{Koskela1} etc.
all the spaces  $  (X, M, \mu, d) $ which satisfied the inequalities (1.1) or (1.2) for each functions
$ \{ u \} $ having the weak gradient correspondingly as a {\it Poincar\'e-Lebesgue spaces} or
{\it Poincar\'e-Lipshitz spaces}.  \par

\vspace{3mm}

 As for the constants $ s, K_P(s,p), K_L(s,p).  $ The value $  s  $ may be defined from the following
condition  (if there exists)

$$
\inf_{x_1,x_2 \in X} \left\{  \frac{\mu(B(x_1, r_1))}{\mu(B(x_2,r_2))} \right\} \ge C \cdot \left[\frac{r_1}{r_2}\right]^s, \ C = \const > 0,
$$
where as ordinary $ B(x,r) $ denotes a closed ball relative the distance $ d(\cdot, \cdot) $ with the center $ x  $ and radii $ r, \ r > 0: $

$$
B(x,r) = \{ y, \ y \in X, \ d(x,y) \le r \}.
$$
This condition is equivalent to the so-called double condition, see \cite{Heinonen1}, \cite{Koskela1}
and  is closely related with the notion of Ahlfors $Q\ - $ regularity

$$
C_1 r^Q \le \mu(B(x,r)) \le C_2 r^Q, \ C_1, C_2, Q = \const > 0,
$$
see  \cite{Koskela1}, \cite{Laakso1}. \par

 Further, we will understand as a capacity of the values $ K_P(s,p), K_L(s,p) $ its minimal values, namely

$$
 K_P(p,q) \stackrel{def}{=} \sup_{0 <||\nabla u||_p < \infty }
  \left\{ \frac{ \mu(X)^{-1/q} ||u - u_X||_q }{ \diam(X) \ \mu(X)^{-1/p} \ ||\nabla u||_p} \right\}, \eqno(1.3a)
$$

$$
 K_L(s,p) \stackrel{def}{=}
\sup_{0 <||\nabla u||_p < \infty } \left\{ \frac{|u(x) - u(y)|}{ d^{1 - s/p}(x,y) \ \mu(X) \ || \nabla u||_p} \right\}. \eqno(1.3b)
$$

\vspace{3mm}

 We  will denote for simplicity

$$
 s =  \order X = \order (X, M, \mu, d).
$$

\vspace{3mm}

  There  are many publications about grounding of these inequalities under some conditions and about its applications, see,
 for instance, in articles \cite{Chuas1}, \cite{Fazio1}, \cite{Hajlasz1}, \cite{Heinonen1}, \cite{Heinonen2}, \cite{Koskela1},  \cite{Koskela2},
 \cite{Kufner1},  \cite{Ostrovsky101},  \cite{Payne1}, \cite{Rjtva1}, \cite{Wannebo1} and in the classical monographs \cite{Beckenbach1},
 \cite{Hardy1}; see also reference therein.\par

 \vspace{4mm}

{\bf  Our aim is a generalization of the estimation (1.1) and (1.2) on the so-called Bilateral
 Grand Lebesgue Spaces $ BGL = BGL(\psi) = G(\psi), $ i.e. when } $ u(\cdot)
 \in G(\psi) \ $  {\it and to show the precision of obtained estimations.} \par

 \vspace{4mm}

  We recall briefly the definition and needed properties of these spaces.
  More details see in the works \cite{Fiorenza1}, \cite{Fiorenza2}, \cite{Ivaniec1},
   \cite{Ivaniec2}, \cite{Ostrovsky1}, \cite{Ostrovsky2}, \cite{Kozatchenko1},
  \cite{Jawerth1}, \cite{Karadzov1} etc. More about rearrangement invariant spaces
  see in the monographs \cite{Bennet1}, \cite{Krein1}. \par

\vspace{3mm}

For $ b = \const, \ 1 < b \le \infty, $ let $\psi =
\psi(p), p \in [1,b),$ be a continuous positive
function such that there exists a limits (finite or not)
$ \psi(1 + 0) $ and $ \psi(b-0), $  with conditions $ \inf_{p \in (1,b)} \psi(p) > 0 $ and
 $\min\{\psi(1+0), \psi(b-0)\}> 0.$  We will denote the set of all these functions
 as $ \Psi(b) $  and   $  b = \supp \psi. $ \par

The Bilateral Grand Lebesgue Space (in notation BGLS) $  G(\psi; a,b) =
 G(\psi) $ is the space of all measurable
functions $ \ f: R^d \to R \ $ endowed with the norm

$$
||f||G(\psi) \stackrel{def}{=}\sup_{p \in (a,b)}
\left[ \frac{ |f|_p}{\psi(p)} \right], \eqno(1.4)
$$
if it is finite.\par
 In the article \cite{Ostrovsky2} there are many examples of these spaces. \par

  The  $ G(\psi) $ spaces over some measurable space $ (X, M, \mu) $
with condition $ \mu(X) = 1 $  (probabilistic case)
appeared in \cite{Kozatchenko1}.\par
 The BGLS spaces are rearrangement invariant spaces and moreover interpolation spaces
between the spaces $ L_1(R^d) $ and $ L_{\infty}(R^d) $ under real interpolation
method \cite{Astashkin2}, \cite{Carro1}, \cite{Jawerth1}, \cite{Karadzov1}. \par

 It was proved also that in this case each $ G(\psi) $ space coincides only
under some additional conditions: convexity of  the functions $ p \to p \cdot \ln \psi(p), \ b = \infty  $ etc.
\cite{Ostrovsky2} with the so-called {\it exponential Orlicz space,} up to norm equivalence. \par

 In others quoted publications were investigated, for instance,
 their associate spaces, fundamental functions
$\phi(G(\psi; a,b);\delta),$ Fourier and singular operators,
conditions for convergence and compactness, reflexivity and
separability, martingales in these spaces, etc.\par

\vspace{3mm}

{\bf Remark 1.1} If we introduce the {\it discontinuous} function

$$
\psi_r(p) = 1, \ p = r; \psi_r(p) = \infty, \ p \ne r, \ p,r \in (a,b) \eqno(1.5)
$$
and define formally  $ C/\infty = 0, \ C = \const \in R^1, $ then  the norm
in the space $ G(\psi_r) $ coincides formally with the $ L_r $ norm:

$$
||f||G(\psi_r) = |f|_r. \eqno(1.5a)
$$

 Thus, the Bilateral Grand Lebesgue spaces are direct generalization of the
classical exponential Orlicz's spaces and Lebesgue spaces $ L_r. $ \par

\vspace{3mm}

{\bf Remark 1.2.}  Let $ F = \{ f_{\alpha}(x) \}, \ x \in X, \ \alpha \in A  $ be certain family
of numerical functions $  f_{\alpha}(\cdot): x \to R, \ A $ is arbitrary set, such that

$$
\exists b > 1, \  \forall p < b \Rightarrow  \psi_F(p) \stackrel{def}{=} \sup_{\alpha \in A}
||f_{\alpha} (\cdot)||_p  < \infty. \eqno(1.6)
$$
 The function $ p \to \psi_F(p) $ is named ordinary as {\it natural} function for the family  $  F.  $  Evidently,

$$
\forall \alpha \in A \Rightarrow  f_{\alpha}(\cdot) \in G\psi_F
$$
and moreover

$$
\sup_{\alpha \in A}|| f_{\alpha}(\cdot)||G\psi_F = 1. \eqno(1.7)
$$

\vspace{4mm}

The BGLS norm estimates, in particular, Orlicz norm estimates for
measurable functions, e.g., for random variables are used in PDE
\cite{Fiorenza1}, \cite{Ivaniec1}, theory of probability in Banach spaces
\cite{Ledoux1}, \cite{Kozatchenko1},
\cite{Ostrovsky1}, in the modern non-parametrical statistics, for
example, in the so-called regression problem \cite{Ostrovsky1}.\par

 We will denote as ordinary the indicator function
$$
I(A) = I(x \in A) = 1, x \in A, \ I(x \in A) = 0, x \notin A;
$$
 here $ A $ is a measurable set.\par

 Recall, see, e.g. \cite{Bennet1} that the fundamental function $ \phi(\delta,S), \ \delta > 0 $ of arbitrary
rearrangement invariant space $  S  $ over $  (X,M,\mu) $ with norm $ ||\cdot||S $  is

$$
\phi(\delta,S) \stackrel{def}{=} ||I(A) ||S, \ \mu(A) = \delta.
$$
 We have in the case of BGLS spaces

$$
\phi(\delta, G\psi) = \sup_{1 \le p < b} \left[ \frac{\delta^{1/p}}{\psi(p)} \right]. \eqno(1.8)
$$
 This notion play a very important role in the functional analysis, operator theory, theory of interpolation and extrapolation,
theory of Fourier series etc., see again \cite{Bennet1}.  Many properties of the fundamental function for
BGLS spaces with considering of several examples see in the articles  \cite{Ostrovsky2}, \cite{Ostrovsky3}. \par

\vspace{3mm}

{\bf Example 1.1.} Let $  \mu(X) = 1 $  and let

$$
\psi^{(b,\beta)}(p) = (b - p)^{-\beta}, \ 1 \le p < b,
\ b = \const > 1, \ \beta = \const > 0, \eqno(1.9)
$$
then as $ \delta \to 0+  $

$$
\phi(G\psi^{(b,\beta)}, \delta) \sim  (\beta b^2/e)^{\beta} \cdot \delta^{1/b} \cdot |\ln \delta|^{-\beta}. \eqno(1.9a)
$$

\vspace{3mm}

{\bf Example 1.2.} Let again $  \mu(X) = 1 $  and let now

$$
\psi_{(\beta)}(p) = p^{\beta}, \ 1 \le p < \infty,
\ \beta = \const > 0, \eqno(1.10)
$$
then as $ \delta \to 0+  $

$$
\phi(G\psi_{(\beta)}, \delta) \sim \beta^{\beta} |\ln \delta|^{-\beta}. \eqno(1.10a)
$$

\vspace{3mm}

\section{Main result: BGLS estimations for  Poincar\'e-Lebesgue-Riesz version. \\
 The case of probability measure.}

\vspace{3mm}

 We suppose in this section without loss of generality that the measure $ \mu $ is probabilistic:
 $  \mu(X) = 1 $  and that the source tetrad $  (X, M, \mu, d) $  is  Poincar\'e-Lebesgue space. \par
  Assume also that the function $ |\nabla u(x)|, \ x \in X  $  belongs to certain BGLS $ G\psi $ with
 $ \supp \psi = s = \order X > 1;  $  the case when   $ \order \psi = b \ne s  $  may be reduced to considered
 here by transfiguration  $ s':=  \min(b,s). $\par

 The function $ \psi(\cdot)  $ may be constructively introduced as a natural function for one function
$ |\nabla u|:  $

$$
\psi_{(0)}(p) := ||u||W_p^1,
$$
if there exists and is finite for  at least one value $ p $ greatest than one. \par

  Define  the following function from the set $  \Psi $

$$
\nu(q) := \inf_{ p \in (qs/(q+s), s) } \left\{ K_P(p,q) \cdot \psi(p)   \right\}, \ 1 \le q < \infty.\eqno(2.1)
$$

 \vspace{3mm}

{\bf Proposition 2.1.}

$$
 ||u - u_X||G\nu \le  \diam(X) \cdot  ||\nabla u||G\psi, \eqno(2.2)
$$
where the "constant"  $  \diam(X) $ is the best possible. \par

\vspace{3mm}

{\bf Proof.} We can suppose without loss of generality $  ||\nabla u||G\psi = 1, $ then it follows by the direct
definition of the norm in BGLS

$$
||\nabla u||_p \le \psi(p), \ 1 \le p < s. \eqno(2.3)
$$

 The inequality (1.1) may be rewritten in our case as follows:

$$
 ||u - u_X||_q \le K_P(p,q) \ \diam(X) \cdot ||\nabla u||_p,
$$
therefore

$$
 ||u - u_X||_q \le K_P(p,q) \ \diam(X) \cdot \psi(p), \ 1 \le p < s. \eqno(2.4)
$$
 Since the value $  p  $ is arbitrary in the set $ 1 \le p < s, $ we can take the minimum of the
right - hand side of the inequality (2.4):

$$
 ||u - u_X||_q \le \diam(X)  \inf_{1 \le p < s} \left[ K_P(p,q)  \cdot \psi(p) \right] =  \diam(X) \cdot \nu(q),
$$
which is equivalent to the required estimate

$$
||u - u_X||G\nu \le \diam(X)  =  \diam(X)||\nabla u||G\psi.
$$

\vspace{3mm}

The exactness of the constant $  \diam(X)  $ in the inequality (2.2) follows immediately from theorem
2.1 in the article \cite{Ostrovsky103}. \par

\vspace{3mm}

\section{Main result: BGLS estimations for  Poincar\'e-Lebesgue-Riesz version.\\
 The general case of arbitrary measure.}

\vspace{3mm}

 The case when the value $  \mu(X) $  is variable, is more complicated.  As a rule,
in the role of a sets $  X  $ acts balls $  B(x,r), $ see  \cite{Heinonen1}, \cite{Koskela1}. \par

\vspace{3mm}

{\bf Definition 3.1.} We will say that the function $ K_P(p,q), \ 1 \le p < s, 1 \le q < \infty $
{\it allows factorable estimation,} symbolically: $  K_P(\cdot,\cdot) \in AFE, $   iff there exist two functions
$ R = R(p) \in \Psi(s) $ and $ V = V(q) \in G\Psi(\infty)  $ such that

$$
K_P(p,q) \le R(p) \cdot V(q). \eqno(3.1)
$$

\vspace{3mm}

{\bf Theorem 3.1.} Suppose that the source tetrad $  (X, M, \mu, d) $  is again Poincar\'e-Lebesgue space
such that $  K_P(\cdot,\cdot) \in AFE. $   Let $ \zeta = \zeta(q) $ be arbitrary function from the set $ \Psi(\infty). $ \par

 Assume also as before in the second section that the function $ |\nabla u(x)|, \ x \in X  $  belongs to certain
BGLS $ G\psi $ with $ \supp \psi = s = \order X > 1;  $  the case when   $ \order \psi = b \ne s  $  may be reduced to
considered here by transfiguration  $ s':=  \min(b,s). $\par

 Our statement:

$$
\frac{||u - u_X||G(V \cdot \zeta)}{ \phi(G\zeta, \mu(X))} \le \diam(X) \cdot \frac{||\nabla u||G\psi}{\phi(R \cdot \psi, \mu(X))} \eqno(3.2)
$$
and  the "constant"  $  \diam(X) $ in (3.2) is as before the best possible. \par

\vspace{3mm}

{\bf Proof.} Denote and suppose for brevity   $ u^{(0)} = u - u_X, \ \mu = \mu(X),  \ \diam(X) = 1, \ ||\nabla u||G\psi = 1.  $
The last equality imply in particular

$$
||\nabla u||_p \le \psi(p), \ 1 \le p < s. \eqno(3.3)
$$
The inequality (1.1) may be reduced taking into account (3.3) as follows

$$
 \mu^{-1/q} ||u^{(0)}||_q \le R(p) \ V(q) \ \mu(X)^{-1/p} \ \psi(p),
$$
and after dividing by $ \zeta(q) $ and by $ \mu^{-1/q} $

$$
\frac{||u^{(0)}||_q }{V(q) \zeta(q)} \le R(p) \ \psi(p) \ \mu^{-1/p} \ \frac{\mu^{1/q}}{\zeta(q)}. \eqno(3.4)
$$
 We take the supremum from both the sides of (3.4) over $  q  $ using the direct definition of the fundamental
function and norm for BGLS:

$$
||u^{(0)}||G(V \cdot \zeta) \le \frac{R(p) \ \psi(p)}{\mu^{1/p}} \cdot \phi(G \zeta,\mu).\eqno(3.5)
$$
 Since the left-hand side of relation (3.5) does not dependent on the variable $  p, $ we can take the
infinum over  $  p.  $  As long as

$$
\inf_p \left[ \frac{R(p) \ \psi(p)}{\mu^{1/p}} \right] =  \left[ \sup_p \frac{\mu^{1/p}}{R(p) \ \psi(p)} \right]^{-1} =
\left[ \phi(G (R \cdot \psi), \mu) \right]^{-1},
$$
we deduce from (3.5)

$$
\frac{||u^{(0)}||G(V \cdot \zeta)}{\phi(G \zeta,\mu)} \le  \frac{1}{\phi(G (R \cdot \psi), \mu) } = \diam X \cdot
\frac{||\nabla u||G\psi}{\phi(G (R \cdot \psi), \mu) },
$$
Q.E.D.

\vspace{3mm}

\section{Main result: BGLS estimations for  Poincar\'e-Lipshitz version.}

\vspace{3mm}

 Recall that we take the number $  s, \ s > 1  $ to be constant. \par
 We consider in this section the case when $ p \in (s,b), \ s < b = \const \le \infty. $

\vspace{3mm}

{\bf Theorem 4.1.} Suppose the fourth $ (X,M, \mu, d) $ is  Poincar\'e-Lipshitz space and that the function
 $ |\nabla u(x)|, \ x \in X  $  belongs to certain BGLS $ G\psi $ with $ \supp \psi = b. $  Then the function
$  u = u(x) $   satisfies after (possible) redefinition on a set of measure zero the inequality

$$
|u(x) - u(y)| \le \mu(X) \cdot \frac{d(x,y)}{\phi( G(K_L \cdot \psi), d^s(x,y))} \cdot ||\nabla u||G\psi, \eqno(4.1)
$$
or equally

$$
\omega(u,\tau) \le \mu(X) \cdot \frac{\tau}{\phi( G(K_L \cdot \psi),\tau^s) } \cdot ||\nabla u||G\psi, \eqno(4.1a)
$$
and this time the "constant" $ \mu(X) $ is best possible. \par

\vspace{3mm}

{\bf Proof.} Suppose for brevity   $  \ \mu(X) = 1, \ ||\nabla u||G\psi = 1.  $
The last equality imply in particular

$$
||\nabla u||_p \le \psi(p), \ s < p  < b.  \eqno(4.2)
$$

 The function  $ u(\cdot) $ satisfies the inequality (1.2)  after (possible) redefinition
of the function $ u = u(x) $ on a set of measure zero

$$
|u(x) - u(y)| \le K_L(s,p) \  d^{1 - s/p}(x,y) \ \mu(X) \ || \nabla u||_p \le
$$

$$
 K_L(s,p) \psi(p)  \cdot d^{1 - s/p}(x,y). \eqno(4.3)
$$

 The excluding set in (4.3) may be dependent on the value $  p, $ but it sufficient to consider
this inequality only for the rational values $  p $  from the interval $  (s,b). $ \par

 The last inequality may be transformed as follows

$$
\frac{|u(x) - u(y)| }{d(x,y)} \le \frac{K_L(s,p)\cdot \psi(p)}{d^{s/p}(x,y)} = \left[ \frac{d^{s/p}(x,y)}{K_L(s,p)\cdot \psi(p)} \right]^{-1}. \eqno(4.4)
$$
Since the left - hand side of (4.4) does not dependent on the variable $  p, $ we can take the infinum from both all the sides of (4.4):

$$
\frac{|u(x) - u(y)| }{d(x,y)} \le \left[ \sup_p \left\{ \frac{d^{s/p}(x,y)}{K_L(s,p)\cdot \psi(p)} \right\} \right]^{-1} =
$$

$$
\frac{1}{\phi(  G(K_L \cdot \psi), d^s(x,y)) } = \frac{\mu(X)}{\phi(  G(K_L \cdot \psi), d^s(x,y)) },
$$
or equally

$$
|u(x) - u(y)| \le \mu(X) \cdot \frac{d(x,y)}{\phi( G(K_L \cdot \psi), d^s(x,y))} =
$$

$$
\mu(X) \cdot \frac{d(x,y)}{\phi( G(K_L \cdot \psi), d^s(x,y))} \cdot ||\nabla u||G\psi. \eqno(4.5)
$$
 The exactness of the constant $ \mu(X) $ may be proved as before, by mention of the article
\cite{Ostrovsky103}.\par
 This completes the proof of theorem 4.1. \par

 Let us consider two examples. \par

 \vspace{3mm}

 {\bf Example 4.1.}  Suppose in addition to the conditions of  theorem 4.1 that  $ \mu(X) = 1 $ and

$$
 K_L(s,p) \cdot \psi(p) =  \psi^{(b,\beta)}(p) = (b - p)^{-\beta}, \ 1 \le p < b,
\ b = \const > 1, \ \beta = \const > 0.
$$

 We deduce  taking into account the example 1.1 that for almost everywhere values $  (x,y) $ and such that $ d(x,y) \le 1/e  $

$$
|u(x) - u(y)| \le C_1(b,\beta,s) \ d^{1 - 1/b}(x,y) \ |\ln d(x,y)|^{\beta} \cdot ||\nabla u||G\psi.
$$

 \vspace{3mm}

 {\bf Example 4.2.}  Suppose in addition to the conditions of  theorem 4.1 that  $ \mu(X) = 1 $ and

$$
 K_L(s,p) \cdot \psi(p) =  \psi_{(\beta)}(p) = p^{\beta}, \ 1 \le p < \infty,
\ \beta = \const > 0.
$$

 We deduce  taking into account the example 1.2 that  for almost everywhere values $  (x,y) $  and such that $ d(x,y) \le 1/e  $

$$
|u(x) - u(y)| \le C_2(\beta,s) \ d(x,y) \ |\ln d(x,y)|^{\beta} ||\nabla u||G\psi.
$$

 \vspace{3mm}

 \section{Concluding remarks}

 \vspace{3mm}

{\bf A.} It may be interest by our opinion to investigate the {\it weights} generalization of obtained
inequalities, alike as done for the classical Sobolev's case,
see for instance  \cite{Chuas1}, \cite{Kufner1}, \cite{Perez1}, \cite{Rjtva1}. \par

\vspace{3mm}

{\bf B.} The physical  applications of these inequalities, for example, in the uncertainty principle. is described in
the article of C.Fefferman \cite{Fefferman1}. \par


\vspace{4mm}

\begin{thebibliography}{99}

\vspace{4mm}

\bibitem{Adams1}
 D.R.Adams, L.I. Hedberg. {\it Function Spaces and Potential Theory.}
 Springer Verlag, Berlin, Heidelberg, New York, 1996.

\bibitem{Astashkin2}
 {\sc Astashkin S.V.} {\it Some new Extrapolation Estimates for the Scale of
 $ L_p \ - $  Spaces.} Funct. Anal. and Its Appl., v. 37 $ N^o $ 3, (2003),  73-77.

\bibitem{Beckenbach1}
E.F.Beckenbach and R.Bellman. {\it Inequalities.} Kluvner Verlag, (1965),
Berlin-Heidelberg-New York.

\bibitem{Bennet1}
C. Bennet and R. Sharpley, {\it Interpolation of operators.}
Orlando, Academic Press Inc., 1988.

\bibitem{Carro1}
M. Carro and J. Martin, {\it Extrapolation theory for the real
interpolation method.} Collect. Math. {\bf 33}(2002), 163--186.

\bibitem{Chuas1}
S.-K. Chuas. {\it Weighted Sobolev inequalities on domains satisfying the
chain condition.} Proc. Amer. Math. Soc., {\bf 122}(4), (2003), 1181-1190.

\bibitem{Fazio1}
G.Di. Fazio, P.Zambroni. {\it A Fefferman-Poincare type inequality for
Carnot-Caratheodory vector fields.}  Proc. of AMS., v. 130 N 9 (2002), p.
2655-2660.

\bibitem{Fefferman1}
C.Fefferman. {\it The uncertainty principle.} Bull. Amer. Math. Soc.,
 {\bf 9,} (1983), p. 129-206.

\bibitem{Fiorenza1}
A. Fiorenza. {\it Duality and reflexivity in grand Lebesgue
spaces.} Collect. Math. {\bf 51}(2000), 131-148.

\bibitem{Fiorenza2}
A. Fiorenza and G.E. Karadzhov, {\it Grand and small Lebesgue
spaces and their analogs.} Consiglio Nationale Delle Ricerche,
Instituto per le Applicazioni del Calcoto Mauro Picone", Sezione
di Napoli, Rapporto tecnico 272/03(2005).

\bibitem{Hajlasz1}
{\sc P. Hajlasz and P. Koskela.} {\it Sobolev meets Poincare.} C. R. Acad. Sci. Paris, 320, (1995),

\bibitem{Hardy1}
G.H.Hardy, J.E. Litlewood and G.P\'olya. {\it Inequalities.}
Cambridge, University Press,  (1952).

\bibitem{Heinonen1}
{\sc J. Heinonen and P. Koskela.} {\it Quasiconformal maps on metric spaces with controlled
geometry.} Acta Math. 181 (1998), 1-61.

\bibitem{Heinonen2}
{\sc J. Heinonen and P. Koskela.} {\it A note on Lipschitz functions, upper gradients and the
Poincar\'e inequality.} New Zealand J. Math. 28 (1999), 37-42.

\bibitem{Ivaniec1}
T. Iwaniec and C. Sbordone, {\it On the integrability of the
Jacobian under minimal hypotheses.} Arch. Rat.Mech. Anal., {\bf 119}, (1992), 129-143.

\bibitem{Ivaniec2}
T. Iwaniec, P. Koskela and J. Onninen, {\it Mapping of Finite
Distortion: Monotonicity and Continuity.} Invent. Math. {\bf 144}, (2001), 507-531.

\bibitem{Jawerth1}
B. Jawerth and M. Milman, {\it Extrapolation theory with
applications.} Mem. Amer. Math. Soc. {\bf 440,} (1991).

\bibitem{Karadzov1}
G.E. Karadzhov and M. Milman, {\it Extrapolation theory: new
results and applications.} J. Approx. Theory, {\bf 113,} (2005), 38-99.

\bibitem{Koskela1}
{\sc P. Koskela. } {\it Upper gradients and Poincare inequalities.} Lectures in Trento in June, 1999.

\bibitem{Koskela2}
{\sc Pekka Koskela.} {\it Sobolev Spaces and Quasiconformal Mappings on Metric Spaces.}
 Internet publications, 2014, PDF.

\bibitem{Kozatchenko1}
Yu.V. Kozatchenko and E.I. Ostrovsky, {\it Banach spaces of random
variables of subgaussian type.} Theory Probab. Math. Stat., Kiev,
1985,  42-56 (Russian).

\bibitem{Krein1}
S.G. Krein, Yu. V. Petunin and E.M. Semenov, {\it Interpolation of
Linear operators.} New York, AMS, 1982.

\bibitem{Kufner1}
A.Kufner. {\it Weighted Sobolev Spaces.} John Wiley Dons, 1985.

\bibitem{Ledoux1}
M. Ledoux and M. Talagrand. {\it Probability in
 Banach Spaces.} Springer, Berlin, 1991.

\bibitem{Mitrinovich1}
{\sc D.S.Mitrinovi\'c, J.E. Pe\^cari\'c and A.M.Fink.} {\it Inequalities Involving
 Functiona and Their Integrals and Derivatives.} Kluvner Academic Verlag,
 (1996), Dorderecht, Boston, London.

\bibitem{Laakso1}
{\sc T. Laakso.}  {\it Ahlfors Q-regular spaces with arbitrary Q admitting weak Poincare inequality.}
Geom. Funct. Anal. 10 (2000), 111-123.

\bibitem{Ostrovsky3}
{\sc E. Liflyand, E. Ostrovsky and L. Sirota.} {\it Structural properties of Bilateral Grand Lebesque Spaces.}
Turk. Journal of Math., {\bf 34,} (2010),  207-219.

\bibitem{Ostrovsky1}
E.I. Ostrovsky, {\it Exponential Estimations for Random Fields.}
Moscow-Obninsk, OINPE, 1999 (Russian).

\bibitem{Ostrovsky2}
E. Ostrovsky and L.Sirota. {\it Moment Banach spaces: theory and applications.}
HAIT Journal of Science and Engeneering, {\bf C}, Volume 4, Issues 1-2,
pp. 233-262, (2007).

\bibitem{Ostrovsky101}
{\sc E. Ostrovsky and L.Sirota.} {\it Poincar\'e inequalities for Bilateral Grand Lebesgue Spaces.}
arXiv:0908.0546v1 [math.FA] 4 Aug 2009

\bibitem{Ostrovsky102}
{\sc E. Ostrovsky and L.Sirota.} {\it Module of continuity for the functions
belonging to the Sobolev-Grand Lebesgue Spaces.}
arXiv:1006.4177v1 [math.FA] 21 Jun 2010

\bibitem{Ostrovsky103}
{\sc E. Ostrovsky and L.Sirota.} {\it Boundedness of operators in Bilateral Grand Lebesgue spaces,
with exact and weakly exact constant calculation.  }
arXiv:1104.2963v1 [math.FA] 15 Apr 2011

\bibitem{Ostrovsky104}
{\sc Ostrovsky E., Rogover E. and Sirota L.} {\it Wirtinger-type inequalities for some rearrangement invariant spaces.}
arXiv:1001.5279v1 [math.FA] 28 Jan 2010

\bibitem{Payne1}
{\sc Payne L.E., Weiberger H.F. } (1960) {\it An optimal Poincar\'e inequality for convex domain.}
Archive for Rational Mechanics and Analysis. V.3, 286-292.

\bibitem{Perez1}
{\sc C.Perez (joint work  with A.Lerner, S.Ombrosi, K.Moen and  R.Torres).} {\it   Sharp
Weighted Bound for Zygmund Singular Integral Operators and Sobolev Inequalities.}
In: "Oberwolfach Reports", Vol. Nunber 3, p. 1828 - 1830; EMS Publishing House,
ETH - Zentrum FLIC1, CH-8092, Zurich, Switzerland.

\bibitem{Poincare1}
{\sc Poincar\'e H.} {\it Sur le \'equations de la physique math\'ematique.} Rend. Circ.
Mat. Palermo 8 (1894), 57-156.; or "Ooevres de Henry Poincar\'e", Paris, (1954),
p. 123-196.

\bibitem{Rjtva1}
Rjtva Hurri Syr\"a. {\it A weighted Poincare inequality with a doubling weight}.
Proc. of the AMS, v. 126 N 2, (1998), p. 542-546.

\bibitem{Takahasi1}
{\sc  S.-E. Takahasi and M. Takeshi.} {\it A note on Wirtinger-Beesack's integral inequalities.}
Nonlinear analysis and convex analysis. Math. Inequal. Appl., {\bf 6}, (2003), no 2,
pp. 277-282.

\bibitem{Takahasi2}
{\sc M. Tsukada, T. Miura, S.Wada, Y.Takahashi and S.E. Takahasi.} {\it On Wirtinger type integral inequalities. }
Nonlinear analysis and convex convex analysis. Yokogama Publ., Yokogama, (2004), v. 5, pp. 541-549.

 \bibitem{Wannebo1}
{\sc A.Wannebo.} {\it Hardy inequalities and imbeddings in domains generalizing
 $ C^{0,\lambda} $ domains.} Proc. Amer. Math. Soc., {\bf 117}, (1993), 449-457;
 {\bf MR} 93d:46050.

\end{thebibliography}
\end{document}